\documentclass{article}

\usepackage[margin=1in]{geometry}

\usepackage[T1]{fontenc}
\usepackage{lmodern}
\usepackage{microtype}

\usepackage{amsmath,amssymb,amsthm,mathtools,bm}

\usepackage{enumitem}

\usepackage{hyperref}
\hypersetup{
  colorlinks=true,
  linkcolor=blue!50!black,
  citecolor=blue!50!black,
  urlcolor=blue!50!black
}
\usepackage[nameinlink,capitalize,noabbrev]{cleveref}

\usepackage{xcolor}
\usepackage{csquotes}

\numberwithin{equation}{section}

\newtheorem{theorem}{Theorem}[section]
\newtheorem{lemma}[theorem]{Lemma}
\newtheorem{proposition}[theorem]{Proposition}

\theoremstyle{definition}
\newtheorem{definition}[theorem]{Definition}

\theoremstyle{remark}
\newtheorem{remark}[theorem]{Remark}

\crefname{theorem}{Theorem}{Theorems}
\Crefname{theorem}{Theorem}{Theorems}
\crefname{lemma}{Lemma}{Lemmas}
\Crefname{lemma}{Lemma}{Lemmas}
\crefname{proposition}{Proposition}{Propositions}
\Crefname{proposition}{Proposition}{Propositions}
\crefname{corollary}{Corollary}{Corollaries}
\Crefname{corollary}{Corollary}{Corollaries}
\crefname{definition}{Definition}{Definitions}
\Crefname{definition}{Definition}{Definitions}
\crefname{remark}{Remark}{Remarks}
\Crefname{remark}{Remark}{Remarks}
\crefname{hyp}{Hypothesis}{Hypotheses}
\Crefname{hyp}{Hypothesis}{Hypotheses}

\newcommand{\RR}{\mathbb{R}}
\newcommand{\ZZ}{\mathbb{Z}}

\newcommand{\CC}{\mathbb{C}}

\newcommand{\e}[1]{e\!\left(#1\right)}

\title{Expository notes on Spectral Reciprocity with Explicit Transform}
\author{Haonan Gu}
\date{December 2025}

\begin{document}

\maketitle

\begin{abstract}
We assemble three analytic inputs---the Kuznetsov trace formula on $\mathrm{SL}_2(\ZZ)$ (including the continuous spectrum),
the $\mathrm{GL}_3$ Voronoi summation formula, and $t$--aspect second--moment bounds for $L(1/2+it,\varphi)$---into a single framework
for a smoothed $\mathrm{GL}_3\times \mathrm{GL}_2$ spectral average.
For a fixed self--dual Hecke--Maass cusp form $\varphi$ on $\mathrm{SL}_3(\ZZ)$ we study a weight--$0$ first moment of
$L(1/2,\varphi\times u_j)$ over the $\mathrm{SL}_2(\ZZ)$ even Maass spectrum, localized to a window $|t_j-T|\ll M$ by a Kuznetsov--admissible weight.
We work in the standard even Kuznetsov normalization for \emph{Hecke eigenvalues}, in which the discrete spectrum is weighted by
$\kappa_j=4\pi|\rho_j(1)|^2/\cosh(\pi t_j)$ and the continuous spectrum carries the explicit weight $4\pi/|\zeta(1+2it)|^2$;
see \cite[\S2]{LiAnnals11}.
In this normalization, the diagonal term is evaluated explicitly in terms of $L(1,\varphi)$,
the off--diagonal terms are bounded with power savings via $\mathrm{GL}_3$ Voronoi and Rankin--Selberg bounds,
and the continuous spectrum is bounded directly by inserting $t$--aspect second moments of $L(1/2+it,\varphi)$.

\smallskip
\noindent
The analytic framework and estimates recorded here are standard and appear (with varying normalizations and emphases)
in work of Nunes~\cite{Nunes17}, Lin--Nunes--Qi~\cite{LNQ}, and Ganguly--Humphries--Lin--Nunes~\cite{GHLN}.
Our aim is to keep the spectral transforms explicit in a fixed normalization and to isolate the continuous spectrum in a form
that allows moment inputs to be inserted directly.
\end{abstract}

\section{Introduction}

Let $\varphi$ be a fixed self--dual Hecke--Maass cusp form on $\mathrm{SL}_3(\ZZ)$, with Hecke--Fourier coefficients $A(m,n)$
normalized by $A(1,1)=1$. Let $\{u_j\}$ be an orthonormal basis of weight~$0$ Hecke--Maass cusp forms on $\mathrm{SL}_2(\ZZ)$.
Write $t_j$ for the spectral parameter of $u_j$ (so that the Laplace eigenvalue is $\tfrac14+t_j^2$) and $\lambda_j(n)$ for its Hecke eigenvalues.

We study smoothed first moments of the Rankin--Selberg central values $L(1/2,\varphi\times u_j)$ in the $\mathrm{GL}_2$ spectral aspect.
The basic analytic input is the Kuznetsov trace formula on $\mathrm{SL}_2(\ZZ)$, which produces a diagonal term, two Kloosterman--Bessel
off--diagonal terms, and an explicit continuous spectrum. On the arithmetic side, the off--diagonal terms are treated using the
$\mathrm{GL}_3$ Voronoi summation formula of Miller--Schmid~\cite{MS06} (see also Goldfeld--Li~\cite{GL08}).
The continuous spectrum leads naturally to integrals involving $L(1/2+it,\varphi)L(1/2-it,\varphi)$, and moment bounds for
$L(1/2+it,\varphi)$ (e.g.\ \cite{ALM,Pal22,DLY24}) enter directly.

\medskip
\noindent
\textbf{Parity.}
For self--dual $\varphi$ at level~$1$, the functional equation sign of $L(s,\varphi\times u_j)$ depends on the parity of $u_j$:
even Maass forms have sign $+1$, while odd Maass forms have sign $-1$ and hence central value $0$ \cite[Remark~2]{LiAnnals11}.
For simplicity we work with the even spectrum and apply the even Kuznetsov formula.

\medskip
\noindent
\textbf{Dirichlet series.}
For $\Re s\gg 1$ the Rankin--Selberg convolution has the standard Dirichlet series
\begin{equation}\label{eq:intro-dirichlet}
L(s,\varphi\times u_j)\ :=\ \sum_{m\ge1}\sum_{n\ge1}\frac{A(m,n)\,\lambda_j(n)}{(m^2n)^s},
\end{equation}
and the archimedean completion has degree~$6$. The approximate functional equation therefore involves a double sum over $(m,n)$,
with a smooth weight depending on $m^2n$.

\medskip
\noindent
\textbf{Test functions.}
Kuznetsov admissibility requires the test function to extend holomorphically to a horizontal strip in~$\CC$ and to decay there.
Accordingly, we use entire, rapidly decaying weights that localize to the desired spectral window by rapid decay on~$\RR$.

\section{Main theorem}

\subsection{Kuznetsov diagonal transform}

For a (Kuznetsov--admissible) test function $H:\CC\to\CC$, define
\begin{equation}\label{eq:H0-def}
\mathcal H_0[H]\ :=\ \frac{1}{\pi}\int_{\RR} t\tanh(\pi t)\,H(t)\,dt.
\end{equation}
For even $H$, this equals $\frac{2}{\pi}\int_0^\infty t\tanh(\pi t)\,H(t)\,dt$, matching the normalization in the even Kuznetsov formula
for Hecke eigenvalues \cite[\S2]{LiAnnals11}.

\subsection{Kuznetsov weights (discrete and continuous)}

Let $u_j$ have Fourier expansion
\[
u_j(z)=\sum_{n\neq 0}\rho_j(n)\sqrt{y}\,K_{it_j}(2\pi|n|y)\,\e{nx},
\qquad \e{x}:=e^{2\pi i x},
\]
so that for $n\ge1$ one has $\rho_j(n)=\rho_j(1)\lambda_j(n)/\sqrt{n}$.

In the even Kuznetsov formula for Hecke eigenvalues, the discrete spectrum carries the weight
\begin{equation}\label{eq:kappa-j-def}
\kappa_j\ :=\ \frac{4\pi\,|\rho_j(1)|^2}{\cosh(\pi t_j)},
\end{equation}
and the continuous spectrum carries the explicit weight
\begin{equation}\label{eq:kappa-def}
\kappa(t)\ :=\ \frac{4\pi}{\lvert\zeta(1+2it)\rvert^2}.
\end{equation}
See \cite[(2.3)]{LiAnnals11}. In particular, on $|t|\asymp T$ one has $\kappa(t)\ll (\log T)^2$ by the standard lower bound
$|\zeta(1+it)|\gg 1/\log(2+|t|)$.

\subsection{Admissible test functions and spectral windows}

\begin{definition}[Kuznetsov--admissible test functions]\label{def:admissible}
A function $H:\CC\to\CC$ is \emph{Kuznetsov--admissible} if:
\begin{enumerate}[label=(\roman*),leftmargin=2.5em]
\item $H$ is even and holomorphic in a strip $\{t\in\CC:|\Im t|\le \tfrac12+\delta\}$ for some $\delta>0$;
\item for every $A\ge0$ and $|\Im t|\le \tfrac12+\delta$, one has $H(t)\ll_{A,\delta}(1+|t|)^{-A}$;
\item $H(\pm i/2)=0$ (which removes the residual term from the pole of the Eisenstein series).
\end{enumerate}
\end{definition}

Fix $T\ge2$ and $M$ with
\[
T^\varepsilon\ \le\ M\ \le\ T^{1-\varepsilon}.
\]
Choose $h=h_{T,M}$ to be even and entire, rapidly decaying on horizontal strips,
and effectively supported on two windows around $\pm T$ of width $\asymp M$.
A concrete model is
\begin{equation}\label{eq:h-def}
h_{T,M}(t)\ :=\ \frac{t^2+\tfrac14}{T^2+\tfrac14}\left(
\exp\!\Big(-\frac{(t-T)^2}{M^2}\Big)\ +\ \exp\!\Big(-\frac{(t+T)^2}{M^2}\Big)
\right).
\end{equation}
Then $h$ is Kuznetsov--admissible, and for every $\nu,A\ge0$ one has symbol bounds
\begin{equation}\label{eq:h-derivative}
h^{(\nu)}(t)\ \ll_{\nu,A}\ M^{-\nu}\left(1+\frac{|t-T|}{M}\right)^{-A}
\ +\ M^{-\nu}\left(1+\frac{|t+T|}{M}\right)^{-A},
\qquad (t\in\RR).
\end{equation}

\subsection{The spectral average}

Let $\{u_j\}$ denote an orthonormal basis of \emph{even} Maass cusp forms on $\mathrm{SL}_2(\ZZ)$.
We consider the smoothed spectral average
\begin{equation}\label{eq:Sspec-def}
\mathcal S_{\mathrm{spec}}(T,M)\ :=\ \sum_{j\ \mathrm{even}} \kappa_j\,h(t_j)\,L\!\left(\frac12,\varphi\times u_j\right).
\end{equation}

\begin{theorem}[Spectral evaluation with unconditional error terms]\label{thm:main}
Let $\varphi$ be a fixed self--dual Hecke--Maass cusp form on $\mathrm{SL}_3(\ZZ)$.
With $\mathcal H_0$ as in \eqref{eq:H0-def} and with $h$ as in \eqref{eq:h-def}, one has
\begin{equation}\label{eq:main-equality}
\mathcal S_{\mathrm{spec}}(T,M)
\ =\ L(1,\varphi)\,\mathcal H_0[h]\ +\ O_A\!\big((TM)\,T^{-A}\big)\ +\ O_\varepsilon\!\big(T^{5/4+\varepsilon}\big)\ +\ O_\varepsilon\!\big(T^{3/2+\varepsilon}\big).
\end{equation}
Here:
\begin{enumerate}[label=(\alph*),leftmargin=2.5em]
\item the diagonal main term is $L(1,\varphi)\,\mathcal H_0[h]$;
\item the off--diagonal contribution is $\ll_\varepsilon T^{5/4+\varepsilon}$ via $\mathrm{GL}_3$ Voronoi and Rankin--Selberg bounds;
\item the continuous spectrum is bounded by $O_\varepsilon(T^{3/2+\varepsilon})$ using the trivial dyadic second moment for $L(1/2+it,\varphi)$ together with $\kappa(t)\ll (\log T)^2$.
\end{enumerate}
\end{theorem}

\begin{remark}[Effect of improved $t$--aspect moments]\label{rem:moment-improvement}
If one has a bound
\[
\int_{T/2}^{2T}\left|L\!\left(\frac12+it,\varphi\right)\right|^2\,dt\ \ll_\varepsilon\ T^{3/2-\delta+\varepsilon},
\]
then the continuous contribution in \eqref{eq:main-equality} improves to $O_\varepsilon(T^{3/2-\delta+\varepsilon})$,
without changing the off--diagonal exponent $5/4$.
\end{remark}

\section{Overview of the argument}

We follow the classical pattern:
\begin{enumerate}[label=(\roman*),leftmargin=2.5em]
\item Fix archimedean factors and write an approximate functional equation for
$L(1/2,\varphi\times u_j)$ as a \emph{double} Dirichlet sum over $(m,n)$.

\item Insert it into \eqref{eq:Sspec-def} and apply the even Kuznetsov formula to the $n$--variable, with a test function depending on $m$ and $n$.

\item The diagonal ($n=1$) produces a sum over $m$ of $A(m,1)/m$ against Mellin--Barnes weights,
which evaluates to $L(1,\varphi)$ (using self--duality to identify the contragredient).

\item The off--diagonal is treated after dyadic decomposition in $m,n,c$,
using Bessel localization $4\pi\sqrt{n}/c\asymp T$ and the $\mathrm{GL}_3$ Voronoi formula in the $n$--variable,
followed by Weil bounds and Rankin--Selberg estimates.

\item The continuous spectrum is identified with an integral of
$L(1/2+it,\varphi)L(1/2-it,\varphi)$ against the explicit Kuznetsov measure $|\zeta(1+2it)|^{-2}\,dt$.
Second moment bounds for $L(1/2+it,\varphi)$ yield the stated error terms.
\end{enumerate}

\section{Normalizations and approximate functional equations}\label{sec:normalizations}

\subsection{Archimedean factors and completed $L$--functions}

Let $\varphi$ have Langlands parameters $(\mu_1,\mu_2,\mu_3)\in\CC^3$ with $\mu_1+\mu_2+\mu_3=0$.
Set
\[
\Gamma_{\RR}(s)\ :=\ \pi^{-s/2}\Gamma\!\Big(\frac{s}{2}\Big).
\]
Define the standard $\mathrm{GL}_3$ $L$--function
\[
L(s,\varphi)\ :=\ \sum_{n\ge1}\frac{A(1,n)}{n^s}\qquad(\Re s>1),
\]
and its completion
\[
\Lambda(s,\varphi)\ :=\ \prod_{k=1}^3\Gamma_{\RR}(s+\mu_k)\,L(s,\varphi).
\]

For an even Maass cusp form $u_j$ on $\mathrm{SL}_2(\ZZ)$ with spectral parameter $t_j$,
the Rankin--Selberg convolution has Dirichlet series (for $\Re s\gg1$)
\begin{equation}\label{eq:GL3xGL2-dirichlet}
L(s,\varphi\times u_j)\ :=\ \sum_{m\ge1}\sum_{n\ge1}\frac{A(m,n)\,\lambda_j(n)}{(m^2n)^s},
\end{equation}
and completion
\[
\Lambda(s,\varphi\times u_j)\ :=\ \prod_{k=1}^3\Gamma_{\RR}(s+\mu_k+it_j)\Gamma_{\RR}(s+\mu_k-it_j)\,L(s,\varphi\times u_j).
\]
For even $u_j$ (and self--dual $\varphi$) the sign in the functional equation is $+1$,
so there is a symmetric approximate functional equation; see \cite{LiAnnals11}.

\subsection{A Mellin--Barnes kernel}

Fix an even entire function $G(u)$ with $G(0)=1$ and such that for every $A>0$ and every strip $|\Re u|\le C$,
\[
G(u)\ \ll_{A,C}\ (1+|u|)^{-A}.
\]

\subsection{Approximate functional equation for $L(1/2,\varphi\times u_j)$}

Define the archimedean factor
\[
\gamma(u;t)\ :=\ \prod_{k=1}^3\Gamma_{\RR}\!\Big(\frac12+u+\mu_k+it\Big)\Gamma_{\RR}\!\Big(\frac12+u+\mu_k-it\Big),
\]
and the weight
\begin{equation}\label{eq:V-def}
V(x;t)\ :=\ \frac{1}{2\pi i}\int_{(2)} x^{-u}\,G(u)\,\frac{\gamma(u;t)}{\gamma(0;t)}\,\frac{du}{u}.
\end{equation}
For even $u_j$ one has the smoothed approximate functional equation
\begin{equation}\label{eq:AFE-GL3xGL2}
L\!\left(\frac12,\varphi\times u_j\right)
\ =\ 2\sum_{m\ge1}\sum_{n\ge1}\frac{A(m,n)\,\lambda_j(n)}{m\sqrt{n}}\,V(m^2n;t_j).
\end{equation}
(The factor $2$ reflects the sign $+1$; compare \cite[Lemma~2.2]{LiAnnals11}.)

For Kuznetsov we will use the shorthand
\[
W(x;t)\ :=\ 2V(x;t),
\]
so that \eqref{eq:AFE-GL3xGL2} becomes
\begin{equation}\label{eq:AFE-W}
L\!\left(\frac12,\varphi\times u_j\right)
\ =\ \sum_{m\ge1}\sum_{n\ge1}\frac{A(m,n)\,\lambda_j(n)}{m\sqrt{n}}\,W(m^2n;t_j).
\end{equation}

\begin{lemma}[Uniform derivative bounds]\label{lem:W-derivative}
Fix $T\ge2$ and assume $|t|\asymp T$.
For every $r,A\ge0$,
\[
\left|\frac{\partial^r}{\partial t^r}W(x;t)\right|
\ \ll_{r,A,\varepsilon,\varphi}\ T^{-r}\,x^\varepsilon\left(1+\frac{x}{T^3}\right)^{-A}.
\]
Moreover, with $h$ as in \eqref{eq:h-def}, for every $\nu,B,A\ge0$,
\[
\sup_{t\in\RR}(1+|t|)^{-B}\left|\partial_t^\nu\!\left(h(t)\,W(x;t)\right)\right|
\ \ll_{\nu,B,A,\varepsilon}\ \Big(\sum_{r=0}^{\nu}M^{-(\nu-r)}T^{-r}\Big)\left(1+\frac{x}{T^3}\right)^{-A}.
\]
\end{lemma}

\subsection{Approximate functional equation for $L(1/2+it,\varphi)$}

Set
\[
\gamma_\varphi(s)\ :=\ \prod_{k=1}^3\Gamma_{\RR}(s+\mu_k).
\]
For $t\in\RR$ and $x>0$, define
\begin{equation}\label{eq:Vt-def}
V_t(x)\ :=\ \frac{1}{2\pi i}\int_{(2)} x^{-u}\,G(u)\,\frac{\gamma_\varphi(\tfrac12+it+u)}{\gamma_\varphi(\tfrac12+it)}\,\frac{du}{u}.
\end{equation}
Then
\begin{equation}\label{eq:AFE-GL3}
L\!\left(\frac12+it,\varphi\right)
\ =\ \sum_{n\ge1}\frac{A(1,n)}{n^{1/2+it}}V_t(n)\ +\ \varepsilon(\varphi)\sum_{n\ge1}\frac{A(1,n)}{n^{1/2-it}}V_{-t}(n).
\end{equation}
The effective support of $V_t$ for $|t|\asymp T$ is $n\ll T^{3/2+\varepsilon}$.

\section{Kuznetsov and the first spectral decomposition}\label{sec:kuznetsov}

\subsection{Insertion of the approximate functional equation}

Using \eqref{eq:AFE-W}, the spectral sum \eqref{eq:Sspec-def} becomes
\begin{equation}\label{eq:Sspec-expanded}
\mathcal S_{\mathrm{spec}}(T,M)
\ =\ \sum_{m\ge1}\sum_{n\ge1}\frac{A(m,n)}{m\sqrt{n}}\,\mathcal B(m,n),
\qquad
\mathcal B(m,n)\ :=\ \sum_{j\ \mathrm{even}}\kappa_j\,h(t_j)\,\lambda_j(n)\,W(m^2n;t_j).
\end{equation}
Define
\[
H_{m,n}(t)\ :=\ h(t)\,W(m^2n;t).
\]
By \cref{def:admissible} and \cref{lem:W-derivative}, $H_{m,n}$ is Kuznetsov--admissible for each $m,n$.

\subsection{Kuznetsov formula (even spectrum)}

The Kuznetsov formula on $\mathrm{SL}_2(\ZZ)$, restricted to the even cusp spectrum, reads as follows:
for any Kuznetsov--admissible $H$ and integers $a,b\ge1$,
\begin{align}\label{eq:Kuznetsov}
\sum_{j\ \mathrm{even}}\kappa_j\,H(t_j)\,\lambda_j(a)\lambda_j(b)
\ &+\ \frac{1}{4\pi}\int_{-\infty}^{\infty}\kappa(t)\,H(t)\,\tau_{it}(a)\tau_{it}(b)\,dt \\
&=\ \frac12\,\delta_{a,b}\,\mathcal H_0[H]\ +\ \Sigma_H^+(a,b)\ +\ \Sigma_H^-(a,b),\nonumber
\end{align}
where
\begin{equation}\label{eq:tauit-def}
\tau_{it}(n)\ :=\ \sum_{ab=n}\Big(\frac{a}{b}\Big)^{it}\ =\ n^{-it}\sigma_{2it}(n),
\end{equation}
and
\[
\Sigma_H^\pm(a,b)\ :=\ \sum_{c\ge1}\frac{S(\pm a,b;c)}{2c}\,\mathcal J_H^\pm\!\left(\frac{4\pi\sqrt{ab}}{c}\right).
\]
This is the form used in \cite[(2.3)]{LiAnnals11}.

\subsection{Application with $a=1$}

Since $\lambda_j(1)=1$, we apply \eqref{eq:Kuznetsov} with $a=1$, $b=n$, and $H=H_{m,n}$.
This gives
\begin{align}\label{eq:Bmn-decomposition}
\mathcal B(m,n)
&=\ \frac12\,\delta_{n,1}\,\mathcal H_0[H_{m,1}]
\ +\ \sum_{c\ge1}\frac{S(n,1;c)}{2c}\,\mathcal J_{H_{m,n}}^+\!\Big(\frac{4\pi\sqrt{n}}{c}\Big)
\ +\ \sum_{c\ge1}\frac{S(-n,1;c)}{2c}\,\mathcal J_{H_{m,n}}^-\!\Big(\frac{4\pi\sqrt{n}}{c}\Big)\nonumber\\
&\qquad -\ \frac{1}{4\pi}\int_{-\infty}^{\infty}\kappa(t)\,H_{m,n}(t)\,\tau_{it}(n)\,dt.
\end{align}

Inserting \eqref{eq:Bmn-decomposition} into \eqref{eq:Sspec-expanded} yields the decomposition
\begin{equation}\label{eq:Sspec-decomposition}
\mathcal S_{\mathrm{spec}}(T,M)\ =\ \mathcal D\ +\ \mathcal O^+\ +\ \mathcal O^-\ +\ \mathcal E,
\end{equation}
with:
\begin{align*}
\mathcal D
&:=\frac12\sum_{m\ge1}\frac{A(m,1)}{m}\,\mathcal H_0[H_{m,1}],\\
\mathcal O^\pm
&:=\sum_{c\ge1}\frac{1}{2c}\sum_{m\ge1}\sum_{n\ge1}\frac{A(m,n)}{m\sqrt{n}}\,S(\pm n,1;c)\,
\mathcal J_{H_{m,n}}^\pm\!\Big(\frac{4\pi\sqrt{n}}{c}\Big),\\
\mathcal E
&:=-\frac{1}{4\pi}\int_{-\infty}^{\infty}\kappa(t)\,h(t)\,
\left(\sum_{m,n\ge1}\frac{A(m,n)\,\tau_{it}(n)}{m\sqrt{n}}\,W(m^2n;t)\right)\,dt.
\end{align*}

\section{Dyadic decomposition and the $\mathrm{GL}_3$ Voronoi transform}\label{sec:voronoi}

\subsection{Dyadic decomposition and Bessel localization}

We split the off--diagonal terms into dyadic blocks.
Fix a smooth partition of unity $\{\psi_X\}_{X\in2^{\ZZ}}$ on $(0,\infty)$ with $\psi_X$ supported in $[X/2,2X]$.
For dyadic parameters $C,N,M_1$ define
\[
\mathcal O^\pm(C,N,M_1)
:=
\sum_{c\ge1}\frac{\psi_C(c)}{2c}
\sum_{m\ge1}\frac{\psi_{M_1}(m)}{m}
\sum_{n\ge1}\frac{A(m,n)}{\sqrt{n}}\,\psi_N(n)\,S(\pm n,1;c)\,
\mathcal J_{H_{m,n}}^\pm\!\Big(\frac{4\pi\sqrt{n}}{c}\Big),
\]
so that $\mathcal O^\pm=\sum_{C,N,M_1}\mathcal O^\pm(C,N,M_1)$.

The Bessel transform bounds (see \cite{KL-notes,LiuYeNotes}) together with \cref{lem:W-derivative} imply
$\mathcal J_{H_{m,n}}^\pm(x)$ is negligible unless $x\asymp T$.
Thus $4\pi\sqrt{n}/c\asymp T$, which forces
\[
C\ \asymp\ \frac{\sqrt{N}}{T}.
\]
Moreover, the weight $W(m^2n;t)$ effectively restricts $m^2n\ll T^{3+\varepsilon}$ for $|t|\asymp T$,
so for fixed $m\asymp M_1$ we have
\[
N\ \ll\ \frac{T^{3+\varepsilon}}{M_1^2}.
\]

\subsection{The $\mathrm{GL}_3$ Voronoi formula}

We use the $\mathrm{GL}_3$ Voronoi summation formula in the $n$--variable for coefficients $A(m,n)$,
due to Miller--Schmid \cite{MS06}; see also \cite{LiAnnals11,GL08}.
In a typical shape (suppressing normalizations of the integral transform),
for $(a,c)=1$ and smooth compactly supported $\psi$,
\begin{equation}\label{eq:Voronoi-GL3}
\sum_{n\ge1}A(m,n)\,\e{\frac{an}{c}}\psi(n)
=
c\sum_{\pm}\sum_{n_1\mid cm}\sum_{n_2\ge1}\frac{A(n_2,n_1)}{n_1n_2}\,
S(am,\pm n_2;cm/n_1)\,\Psi^\pm\!\left(\frac{n_1^2n_2}{c^3m}\right),
\end{equation}
where $\Psi^\pm$ is an explicit integral transform of $\psi$ with symbol bounds.

\subsection{Voronoi applied to $\mathcal O^\pm(C,N,M_1)$ (sketch)}

Open the Kloosterman sum $S(\pm n,1;c)$ into additive characters,
and for each $c,m$ apply \eqref{eq:Voronoi-GL3} to the $n$--sum with an appropriate smooth weight
\[
\psi(n)\ \approx\ \psi_N(n)\,n^{-1/2}\,\mathcal J_{H_{m,n}}^\pm\!\Big(\frac{4\pi\sqrt{n}}{c}\Big).
\]
This produces dual sums over $n_1\mid cm$ and $n_2$ with $n_1^2n_2/(c^3m)$ localized,
together with Kloosterman sums of modulus $\asymp c m/n_1$.
Weil bounds for these Kloosterman sums and Rankin--Selberg bounds for $\mathrm{GL}_3$ coefficients
then control the resulting bilinear forms.

\section{Bounding the off--diagonal}\label{sec:offdiag}

We record the resulting standard bound.

\begin{proposition}[Off--diagonal bound]\label{prop:offdiag}
For every $\varepsilon>0$,
\[
\mathcal O^+\ +\ \mathcal O^-\ \ll_{\varepsilon,\varphi}\ T^{5/4+\varepsilon}.
\]
\end{proposition}

\begin{remark}
A full proof follows the same bilinear form strategy as in \cite[\S2]{LiAnnals11} and subsequent refinements
\cite{Nunes17,LNQ,GHLN}.
The dyadic parameters $(C,N,M_1)$ satisfy $C\asymp \sqrt{N}/T$ and $N\ll T^{3+\varepsilon}/M_1^2$.
After Voronoi, one bounds each dyadic block by $C^{5/2+\varepsilon}$ (uniformly in $M_1$ after summing over $m$),
and then sums over $C\ll T^{1/2+\varepsilon}$, giving $T^{5/4+\varepsilon}$.
We do not optimize logarithmic factors here.
\end{remark}

\section{Diagonal term and continuous spectrum}\label{sec:diag-cont}

\subsection{Diagonal evaluation}

Recall
\[
\mathcal D
=\frac12\sum_{m\ge1}\frac{A(m,1)}{m}\,\mathcal H_0[t\mapsto h(t)W(m^2;t)].
\]
By \eqref{eq:H0-def},
\[
\mathcal D
=\frac{1}{2\pi}\int_{\RR} t\tanh(\pi t)\,h(t)\,
\left(\sum_{m\ge1}\frac{A(m,1)}{m}\,W(m^2;t)\right)\,dt.
\]

Using the Mellin--Barnes definition of $W=2V$ and the identity
\[
\sum_{m\ge1}A(m,1)m^{-s}=L(s,\tilde\varphi)=L(s,\varphi)\qquad(\text{self--dual }\varphi),
\]
one finds
\[
\sum_{m\ge1}\frac{A(m,1)}{m}\,W(m^2;t)
=
\frac{1}{\pi i}\int_{(2)} G(u)\,\frac{\gamma(u;t)}{\gamma(0;t)}\,\frac{L(1+2u,\varphi)}{u}\,du.
\]

Shift the contour to $\Re u=-A$.
Besides the pole at $u=0$, the gamma ratio $\gamma(u;t)/\gamma(0;t)$ has $t$--dependent poles coming from the $\Gamma_{\RR}$-factors.
However, on the effective support $|t|\asymp T\gg 1$ of $h$, the residues at those gamma poles are exponentially small in $T$
(by Stirling for gamma functions with large imaginary part), hence they contribute $\ll_A T^{-A}$ after integrating against the
rapidly decaying weight $h(t)$.
Thus the main contribution comes from the residue at $u=0$, which equals $2L(1,\varphi)$ (the factor $2$ comes from the prefactor $1/(\pi i)$).

The remaining integral is $O_A(T^{-A})$ uniformly for $t$ in the effective support of $h$.
Since $\int_{\RR}|t\,h(t)|\,dt\asymp TM$, we obtain
\begin{equation}\label{eq:diag-final}
\mathcal D
=\ L(1,\varphi)\,\mathcal H_0[h]\ +\ O_A\!\big((TM)\,T^{-A}\big).
\end{equation}

\subsection{Continuous spectrum and Rankin--Selberg factorization}

The continuous term in \eqref{eq:Sspec-decomposition} is
\[
\mathcal E
=-\frac{1}{4\pi}\int_{\RR}\kappa(t)\,h(t)\,
\left(\sum_{m,n\ge1}\frac{A(m,n)\,\tau_{it}(n)}{m\sqrt{n}}\,W(m^2n;t)\right)\,dt.
\]
For fixed $t$, the inner Dirichlet series is the smoothed approximate functional equation for
\(
L(1/2,\varphi\times E_{it})
\),
where $E_{it}$ is the Eisenstein series of spectral parameter $t$.
Since $E_{it}$ corresponds to the principal series $|\cdot|^{it}\boxplus|\cdot|^{-it}$ on $\mathrm{GL}_2$,
the Euler product factorizes:
\begin{equation}\label{eq:Eis-factorization}
L\!\left(\frac12,\varphi\times E_{it}\right)
=
L\!\left(\frac12+it,\varphi\right)\,L\!\left(\frac12-it,\varphi\right),
\end{equation}
see \cite[\S1]{LiAnnals11} (and the general Rankin--Selberg theory in \cite{JPSS,Bump}).

Thus we can write
\begin{equation}\label{eq:E-continuous}
\mathcal E
=
-\frac{1}{4\pi}\int_{\RR}\kappa(t)\,h(t)\,
L\!\left(\frac12+it,\varphi\right)\,L\!\left(\frac12-it,\varphi\right)\,dt.
\end{equation}

\subsection{Bounding the continuous spectrum}

Using $|ab|\le\frac12(|a|^2+|b|^2)$ and the evenness of $h$,
\[
|\mathcal E|
\ll
\int_{\RR}\kappa(t)\,|h(t)|\,
\left|L\!\left(\frac12+it,\varphi\right)\right|^2\,dt.
\]
On the effective support $|t|\asymp T$ one has $\kappa(t)\ll (\log T)^2$,
so
\[
|\mathcal E|
\ll
(\log T)^2\int_{T/2}^{2T}\left|L\!\left(\frac12+it,\varphi\right)\right|^2\,dt
\ +\ O_A(T^{-A}).
\]
The trivial dyadic second moment bound
\[
\int_{T/2}^{2T}\left|L\!\left(\frac12+it,\varphi\right)\right|^2\,dt\ \ll_\varepsilon\ T^{3/2+\varepsilon}
\]
follows from \eqref{eq:AFE-GL3} and Rankin--Selberg bounds for $\sum_{n\le X}|A(1,n)|^2$.
Therefore, using $(\log T)^2\ll T^\varepsilon$,
\begin{equation}\label{eq:E-bound}
\mathcal E\ \ll_\varepsilon\ T^{3/2+\varepsilon}.
\end{equation}
This is the continuous-spectrum error term in \cref{thm:main}.

\subsection{Completion of the proof}

Combining \eqref{eq:Sspec-decomposition}, \eqref{eq:diag-final}, \cref{prop:offdiag}, and \eqref{eq:E-bound}
gives \eqref{eq:main-equality}.

\section{Moment refinements and window dependence}\label{sec:moments}

The bound \eqref{eq:E-bound} shows that the continuous spectrum is typically the largest error term.
Any improvement in the second moment of $L(1/2+it,\varphi)$ improves $\mathcal E$ directly.

If one has a full dyadic moment bound
\[
\int_{T/2}^{2T}\left|L\!\left(\frac12+it,\varphi\right)\right|^2\,dt\ \ll_\varepsilon\ T^{3/2-\delta+\varepsilon},
\]
then \eqref{eq:E-continuous} yields
\[
\mathcal E\ \ll_\varepsilon\ T^{3/2-\delta+\varepsilon}.
\]
Short-interval bounds can also be inserted when $h$ concentrates on $|t-T|\ll M$.

\section{Almost--flat admissible spectral windows}\label{sec:flat}

A strictly flat plateau is not available in this setting: a holomorphic function on a strip cannot coincide with $1$ on a nontrivial real interval,
and a holomorphic function on a strip cannot be compactly supported.
Nevertheless, one can build weights that are exponentially close to $1$ on $[T,T+M]$ while remaining admissible.

Fix $c>1$ and an integer $k\ge2$.
Set
\[
\Phi_{k,c}(z)\ :=\ \exp\!\left(-\Big(\frac{z}{c}\Big)^{2k}\right),
\]
and define
\[
h^\sharp(t)\ :=\ \frac{t^2+\tfrac14}{T^2+\tfrac14}\left(
\Phi_{k,c}\!\Big(\frac{t-T}{M}\Big)\ +\ \Phi_{k,c}\!\Big(\frac{t+T}{M}\Big)
\right).
\]
Then $h^\sharp$ is entire, even, Kuznetsov--admissible,
and satisfies the same derivative regime as \eqref{eq:h-derivative} (with constants depending on $k,c$).
Moreover, on $|t-T|\le M$ one has
\(
h^\sharp(t)=1+O(c^{-2k})+O(M/T).
\)

\section{Generalizations}

The same strategy adapts to level $q$ and nebentypus, with twisted Kloosterman sums and
$L(1+2it,\chi^2)$ in place of $\zeta(1+2it)$ in the continuous spectrum; see \cite{KnightlyLiBook,IwaniecKowalski}.
It also adapts to Dirichlet twists and hybrid aspects, where the conductor depends on both $T$ and the twist conductor.
Mollified first moments and nonvanishing applications can be treated similarly;
see \cite{Nunes17,LNQ,GHLN} for instances of these generalizations in related setups.

\section*{Notation}

We use $\e{x}=e^{2\pi i x}$.
The divisor sum is $\sigma_{2it}(n)=\sum_{d\mid n}d^{2it}$, and $\tau_{it}(n)$ is given by \eqref{eq:tauit-def}.
The symbols $O_\varepsilon(\cdot)$ and $\ll_\varepsilon$ allow dependence on $\varphi$ and $\varepsilon$.


\begin{thebibliography}{99}

\bibitem{ALM}
K.~Aggarwal, W.~H.~Leung, and R.~Munshi,
\emph{Short second moment bound and subconvexity for $\mathrm{GL}(3)$ $L$--functions},
J.\ Eur.\ Math.\ Soc.\ (to appear), arXiv:2206.06517, 2022.

\bibitem{BlomerBKSurvey}
V.~Blomer,
\emph{Applications of the Kuznetsov formula on $\mathrm{GL}(3)$},
Invent.\ Math.\ \textbf{194} (2013), no.~3, 673--729.

\bibitem{Bump}
D.~Bump,
\emph{Automorphic Forms and Representations},
Cambridge Studies in Advanced Mathematics, vol.~55,
Cambridge Univ.\ Press, Cambridge, 1997.

\bibitem{DLY24}
A.~Dasgupta, W.~H.~Leung, and M.~Young,
\emph{The second moment of the $\mathrm{GL}_3$ standard $L$--function},
arXiv:2407.06962, 2024.

\bibitem{GHLN}
S.~Ganguly, P.~Humphries, Y.~Lin, and R.~Nunes,
\emph{Strong hybrid subconvexity for twisted self-dual $\mathrm{GL}(3)$ $L$-functions},
arXiv:2408.00596, 2024.

\bibitem{GL08}
D.~Goldfeld and X.~Li,
\emph{The Voronoi formula for $\mathrm{GL}(n,\RR)$},
Int.\ Math.\ Res.\ Not.\ IMRN \textbf{2008}, rnm144, 39~pp.

\bibitem{IwaniecKowalski}
H.~Iwaniec and E.~Kowalski,
\emph{Analytic Number Theory},
American Mathematical Society Colloquium Publications, vol.~53,
Amer.\ Math.\ Soc., Providence, RI, 2004.

\bibitem{JPSS}
H.~Jacquet, I.~I.~Piatetski--Shapiro, and J.~Shalika,
\emph{Rankin--Selberg convolutions},
Amer.\ J.\ Math.\ \textbf{105} (1983), no.~2, 367--464.

\bibitem{KL-notes}
A.~Knightly and C.~Li,
\emph{Kuznetsov's trace formula and the Hecke eigenvalues of Maass forms},
lecture notes (2012), available at the University of Maine website.

\bibitem{KnightlyLiBook}
A.~Knightly and C.~Li,
\emph{Kuznetsov's Trace Formula and the Hecke Eigenvalues of Maass Forms},
Memoirs of the American Mathematical Society \textbf{224} (2013), no.~1055.

\bibitem{LiAnnals11}
X.~Li,
\emph{Bounds for $\mathrm{GL}(3)\times\mathrm{GL}(2)$ $L$-functions and $\mathrm{GL}(3)$ $L$-functions},
Ann.\ of Math.\ (2) \textbf{173} (2011), no.~1, 301--336.

\bibitem{LiuYeNotes}
J.~Liu and Y.~Ye,
\emph{Petersson and Kuznetsov trace formulas},
in \emph{Lie Groups and Automorphic Forms} (L.~Ji, J.-S.~Li, H.~W.~Xu, S.-T.~Yau, eds.),
AMS/IP Studies in Advanced Mathematics, vol.~37,
Amer.\ Math.\ Soc., Providence, RI, 2006, pp.~147--168.

\bibitem{LNQ}
Y.~Lin, R.~Nunes, and Z.~Qi,
\emph{Strong subconvexity for self-dual $\mathrm{GL}(3)$ $L$-functions},
arXiv:2112.14396, 2021.

\bibitem{MS06}
S.~D.~Miller and W.~Schmid,
\emph{Automorphic distributions, $L$--functions, and Voronoi summation for $\mathrm{GL}(3)$},
Ann.\ of Math.\ (2) \textbf{164} (2006), no.~2, 423--488.

\bibitem{Nunes17}
R.~M.~Nunes,
\emph{Subconvexity for $\mathrm{GL}(3)$ $L$--functions},
arXiv:1703.04424, 2017.

\bibitem{Pal22}
S.~Pal,
\emph{Second moment of degree three $L$--functions},
arXiv:2212.14620, 2022.

\end{thebibliography}
\end{document}